\documentclass[12 pt,a4paper]{article}

\usepackage{amsfonts,amsmath,amsthm,amssymb,setspace}
\newtheorem{Proposizione}[equation]{Proposition}
\newtheorem{Definizione}[equation]{Definition}
\newtheorem{Teorema}[equation]{Theorem}
\newtheorem{Remark}{Remark}
\newtheorem{Corollario}[equation]{Corollary}
\numberwithin{equation}{section}
\begin{document}


\author{Luciano Mari \and Marco Rigoli}
\title{\textbf{Maps from Riemannian manifolds into non-degenerate Euclidean
cones}}

\date{}
\maketitle
\scriptsize \begin{center} Dipartimento di Matematica,
Universit\`a
degli studi di Milano,\\
Via Saldini 50, I-20133 Milano (Italy)\\
E-mail addresses: luciano.mari@unimi.it, marco.rigoli@unimi.it
\end{center}

\begin{abstract}\footnote{\textbf{Mathematic subject classification 2000}: primary
53C42, 35B50; secondary 53C21.\par
\textbf{ \ Keywords}: maximum principles, harmonic maps, isometric
immersion, Riemannian manifold.}
Let $M$ be a connected, non-compact $m$-dimensional Riemannian
manifold. In this paper we consider smooth maps $\varphi: M
\rightarrow \mathbb{R}^n$ with images inside a non-degenerate
cone. Under quite general assumptions on $M$, we provide a lower
bound for the width of the cone in terms of the energy and the
tension of $\varphi$ and a metric parameter. As a side product, we
recover some well known results concerning harmonic maps, minimal
immersions and K\"ahler submanifolds. In case $\varphi$ is an
isometric immersion, we also show that, if $M$ is sufficiently
well-behaved and has non-positive sectional curvature,
$\varphi(M)$ cannot be contained into a non-degenerate cone of
$\mathbb{R}^{2m-1}$.
\end{abstract}

\normalsize
\section{Introduction}
In a famous paper Omori, \cite{O}, studied minimal immersions of
manifolds into cones of Euclidean space and in the process he
proved (a version of) what is by now known as the Omori-Yau
maximum principle. This powerful tool has since then been applied
to solve several geometric problems from Yau's pioneering work on
the generalized Schwarz lemma, to the study of the group of
conformal diffeomorphisms of a manifold, \cite{PRS2}, to the
theory of submanifolds (see also the very recent papers \cite{APD}
and \cite{AG}) and so on; for a short account one can consult
\cite{PRS1} and the references therein. The Omori-Yau principle
can be briefly stated in the following form:
\begin{quote}
Let $(M,\langle,\rangle)$ be a complete Riemannian manifold and
let $u\in C^2(M)$ be a function satisfying $u^*=\sup_M u
<+\infty$. Then, under appropriate geometric assumptions on the
manifold there exists a sequence $\{x_k\}\subset M$ with the
following properties:
$$
(i) \quad u(x_k)>u^*-\frac{1}{k} \quad , \quad (ii) \quad |\nabla
u|(x_k) <\frac{1}{k} \quad , \quad (iii) \quad \Delta u(x_k)
<\frac{1}{k}
$$
for every $k\in \mathbb{N}$.
\end{quote}
It turns out that in many applications property $(ii)$ plays no
role. This fact suggests to ignore it and to look for possibly
lighter geometric assumptions to ensure the validity of
conclusions $(i)$ and $(iii)$. Indeed, in \cite{PRS2} the authors
show that the sole assumption of \emph{stochastic completeness} of
$(M,\langle,\rangle)$ (that is, the property of the Brownian
motion on $M$ to have infinite intrinsic lifetime) is equivalent
to the simultaneous validity of $(i)$ and $(iii)$ for every $u\in
C^2(M)$, $u^*<+\infty$, on an appropriate sequence.\\
To distinguish between the two situations we shall refer to this
new statement as to the \emph{weak maximum principle}. These two
principles are indeed different. For instance, one immediately
realizes that stochastic completeness does not even imply geodesic
completeness (as an example, consider the standard punctured
Euclidean space $\mathbb{R}^n\backslash \{o\}$).\par
As a matter of fact, geometric applications often lead to consider
differential inequalities of the form
$$
\Delta u \ge b(x)f(u)
$$
for some $b(x)\in C^0(M)$, $f(t)\in C^0(\mathbb{R})$. This and a
number of other considerations have suggested to generalize the
known results on stochastic completeness and the weak maximum
principle to a wider class of elliptic operators. In particular,
to symmetric diffusion operators $L = b^{-1}\mathrm{div}(A\nabla)$
on $L^2(M,bdx)$, where $A,b$ are positive smooth functions on $M$.
Theorem 3.10 of \cite{PRS1} \emph{states the equivalence between
$L$-stochastic completeness} (here $L$ means that the diffusion
having infinite intrinsic lifetime is that generated by $L$)
\emph{and the weak maximum principle for $L$}. For a detailed
account, one can consult \cite{PRS1}, \cite{RSV}, where the
authors deal
even with the case $u^*=+\infty$.\\
Together with diffusions, one can also consider the full Hessian
operator. In this case we state the following definition:
\begin{quote}
The \emph{weak maximum principle for the Hessian operator} holds
on $M$ if for every $u \in C^2(M)$, $u^\ast < +\infty$, there
exists a sequence $\{x_k\} \subset M$ such that
\begin{equation}\label{hessianoprinc}
u(x_k) > u^\ast -\frac{1}{k} \quad ,
\sup_{\scriptsize {\begin{array}{c} v \in T_{x_k} M \\
|v| = 1 \end{array}}} \mathrm{Hess}_{x_k} u (v,v) < \frac{1}{k}
\end{equation}
\end{quote}
It is worth to observe, in view of applications, that
$L$-stochastic completeness and the weak maximum principle for the
Hessian can be granted via the existence of appropriate proper
functions defined in a neighborhood of infinity in $M$. This is
achieved by the use of Khas'minskii test (see \cite{PRS1}, where
for convenience the result is stated for the Laplacian only) and
of its generalization in the case of the Hessian (see again
\cite{PRS1} and also question $33$ in \cite{PRS4}). Observe that
adding the condition $|\nabla u(x_k)| < 1/k$ to
\eqref{hessianoprinc} gives the well known Omori-Yau maximum
principle for the Hessian.\par
In this paper we shall make use of Theorem $1.9$ (and subsequent
remarks), Theorem 3.10 and Proposition $3.18$ of \cite{PRS1} to
prove the geometrical results
that we are now going to describe.\\
\par
Let $(\mathbb{R}^n,\langle,\rangle)$, $n\ge 2$, be the
$n$-dimensional Euclidean space endowed with its canonical flat
metric. Fix an origin $o \in \mathbb{R}^n$ and a unit vector $v\in
\mathbb{S}^{n-1}$. We set $\mathcal{C}_{o,v,\theta}$, shortly
$\mathcal{C}$, to denote the non-degenerate cone with vertex in
$o$, direction $v$ and width $\theta$, $\theta\in (0,\pi/2)$, that
is,
\begin{equation}
\mathcal{C} = \Big\{ z \in \mathbb{R}^n\backslash\{o\} \ : \ \Big\langle \frac{z-o}{|z-o|},v \Big\rangle \ge \cos (\theta)
\Big\}.
\end{equation}
Let $(M,(,))$ be a connected, $m$-dimensional Riemannian manifold
($m\ge 2$), and let
$$
\varphi \ : \quad (M,(,)) \longrightarrow
(\mathbb{R}^n,\langle,\rangle)
$$
be a smooth map. We indicate with $|d\varphi|^2$ the square of the
Hilbert-Schmidt norm of the differential $d\varphi$ (in other
words, twice the energy density of $\varphi$) and with
$\tau(\varphi)$ the tension field of $\varphi$. In case $\varphi$
is an isometric immersion, $|d\varphi|^2 = m$ and $\tau(\varphi)=
mH$, where $H$ is the mean curvature vector. We fix an origin
$q\in M$ and we consider the distance function from $q$,
$r(x)=d(x,q)$. We set $B_r$ for the geodesic ball with radius $r$
centered at $q$.\\
\par
Since the work of Omori, the study of maps with image contained in
a cone has captured the attention of researchers. Indeed,
important progresses in the understanding of their geometry have
been made by a number of authors and in particular, regarding the
content of the present paper, by \cite{BK}, \cite{PRS1} and by
\cite{A}, where through stochastic methods it is proved that if
$M$ is stochastically complete and $\varphi$ is an harmonic map
with energy density bounded from below away from zero,
$\varphi(M)$ cannot be contained into any non-degenerate cone of
$\mathbb{R}^n$. However, all these references limit themselves to
consider the harmonic (or minimal) case.\par
On the contrary, the two main results in this paper aim at
determining a lower bound for the width of the cone containing
$\varphi(M)$ when $\varphi$ is merely a smooth map or an isometric
immersion, where of course we expect the estimate to depend on the
position of
$\varphi(M)$ in the cone.\\
\par
\noindent To state the next theorems, given $\eta
>0$, we define
\begin{equation}\label{Aeta}
A_\eta = \sup_{\scriptsize{ \begin{array}{c} \xi \in (0,1) \\
\alpha \in (0, \min\{1,\eta\sqrt{1-\xi}\}) \end{array}}} \xi \alpha^2 \sqrt{1-\alpha^2}.
\end{equation}
The constant $A_\eta$ can be easily computed, but the actual value
is irrelevant for our purposes. Note also that $A_\eta$ is
non-decreasing as a function of $\eta$.
\begin{Teorema}\label{stimacono}
Let M be a connected, non-compact $m$-dimensional Riemannian
manifold, and let
$$
\varphi \ : \quad (M,(,)) \longrightarrow
(\mathbb{R}^n,\langle,\rangle)
$$
be a map of class $C^2$ such that $|d\varphi(x)|^2>0$ on $M$.
Consider the elliptic operator $L=|d\varphi|^{-2}\Delta$, and
assume that $M$ is $L$-stochastically complete. Let
$\mathcal{C}=\mathcal{C}_{o,v,\theta}$ be a cone with vertex at
$o\in\mathbb{R}^n\backslash\varphi(M)$, let $\pi_v$ be the
hyperplane orthogonal to $v$ passing through $o$ and let
$d(\pi_v,\varphi(M))$ be the Euclidean distance between this
hyperplane and $\varphi(M)$.\par
If $\varphi(M)$ is contained in $\mathcal{C}$, then
\begin{equation}\label{principale}
\cos(\theta) \le \sqrt{\frac{1}{A_1} d(\pi_v,\varphi(M)) \sup_M \Big[\frac{|\tau(\varphi)|}{|d\varphi|^2}\Big]}.
\end{equation}
In case $\varphi$ is an isometric immersion, we can replace $A_1$
with $A_m$ in (\ref{principale}) obtaining a sharper estimate.
\end{Teorema}

\begin{Remark}
\emph{Note that, in case
$$
\sup_M \Big[\frac{|\tau(\varphi)|}{|d\varphi|^2}\Big]=+\infty
$$
and $d(\pi_v,\varphi(M))=0$, that is, $\varphi(M)$ "gathers around the origin $o$", as we shall see in the
proof, we have no restriction on $\theta$.}
\end{Remark}

\begin{Remark}\label{sufficienticond}
\emph{The condition that $L=|d\varphi|^{-2}\Delta$ generates a
conservative diffusion (i.e. $L$-stochastic completeness) is
implied by
$$
\int_0^\xi |d\varphi|^2 dt = +\infty \qquad \mathbb{P}_x \ a.s \text{ for some (hence, every) } x\in M,
$$
where the integral is evaluated along the paths of the Brownian
motion $X_t$ starting from $x$, and $\xi$ is the intrinsic
lifetime of $X_t$. As in the case of $\Delta$-stochastic
completeness, no geodesic completeness is required. If $M$ is
complete, $L$-stochastic completeness can be obtained using
generalized versions of a previous work of Grigor'Yan \cite{G},
see \cite{S} and \cite{PRS1} (Theorem 3.15 and subsequent
propositions). In particular, if there exist $C>0$, $\beta\in
\mathbb{R}$ such that
\begin{equation}\label{casoB}
|d\varphi(x)|^2 \ge \frac{C}{(1+r(x))^\beta} \quad \mathit{on \
}M
\end{equation}
and
\begin{equation}\label{volumeB}
 \frac{r^{1-\beta}}{\log(\mathrm{Vol}(B_r))}
\not \in L^1(+\infty),
\end{equation}
then the weak maximum principle for $L=|d\varphi|^{-2}\Delta$
holds. It is worth to observe that \eqref{volumeB} implies
$\beta\le 2$, but no restriction on nonnegativity of $\beta$ is
needed. In case $\beta=2$, an application of \cite{S} leads to
slightly improving \eqref{volumeB} to
$$
\frac{\log r}{r\log(\mathrm{Vol}(B_r))} \not\in L^1(+\infty).
$$
The converse of the above result is in general false; indeed, it
is well known that, on a radially symmetric model in the sense of
Greene and Wu, \cite{GW}, $\Delta$-stochastic completeness is
equivalent to the condition
\begin{equation}\label{modello}
\frac{\mathrm{Vol}(B_r)}{\mathrm{Vol}(\partial B_r)} \not \in L^1(+\infty),
\end{equation}
which is implied by \eqref{volumeB} with $\beta=0$ but not
equivalent. We remark, in passing, that it is an open problem to
establish whether or not \eqref{modello} is sufficient for
$\Delta$-stochastic completeness on a generic complete manifold.}
\end{Remark}

\begin{Remark}[\textbf{Sharpness of inequality \ref{principale}}]
\emph{Due to the form of \eqref{principale}, we cannot expect the
result to be significant when $\varphi(M)$ is far from $\pi_v$, in
the following sense: for every $M$, $\mathcal{C}$ and $\varphi$
satisfying the assumptions of Theorem \ref{stimacono}, and for
every $k\ge 0$, we can consider the map $\varphi_k = \varphi +
kv$. Then $d(\pi_v,\varphi_k(M))= d(\pi_v,\varphi(M))+k$, while
the other parameters in the RHS of \eqref{principale} remain
unchanged. Therefore, for $k$ sufficiently large inequality
\eqref{principale} becomes meaningless unless $\tau(\varphi)
\equiv 0$. On the contrary, we show with a simple example that,
when $d(\pi_v,\varphi(M))$ is very small, \eqref{principale} is
sharp in the following sense: for every fixed hyperplane $\pi_v$,
and for every origin $o\in \pi_v$, there exists a family of maps
$\varphi_d$, $d>0$ representing $d(\pi_v,\varphi_d(M))$, such
that, if we denote with $\theta_d$ the width of the non-degenerate
tangent cone containing $\varphi_d(M)$,
$$
\frac{\cos^2(\theta_d)}{d} \ge C \qquad \text{when} \quad  d \rightarrow 0^+,
$$
for some constant $C>0$. Indeed, for every fixed $d>0$ consider
the hypersurface $\varphi_d : \mathbb{R}^m \rightarrow
\mathbb{R}^{m+1}$ given by the graph $\varphi(x) = (x, |x|^2 +
d)$, with the induced metric. Indicating with $\pi_v$ the
hyperplane $x_{m+1}=0$, we have by standard calculations
$$
|\tau(\varphi_d)|  = \frac{2m + 8(m-1)|x|^2}{(1+4|x|^2)^{3/2}}
\quad \text{and} \quad |d\varphi_d|^2 = m
$$
Therefore $\sup_M |\tau(\varphi_d)|/|d\varphi_d|^2= 2$. Moreover,
for the tangent cone passing through the origin
$$
\cos^2(\theta_d) = \frac{4d}{1+4d},
$$
thus, since $d \equiv d(\pi_v,\varphi_d(M))$, we reach the desired
conclusion. }\end{Remark}
As an immediate consequence of Theorem \ref{stimacono}, we recover
Atsuji's result \cite{A}:
\begin{Corollario}\label{atsuji}
Let $\varphi : M \rightarrow \mathbb{R}^n$ be harmonic and such
that $|d\varphi|^2 \ge C$ for some positive constant $C$. If $M$
is stochastically complete, then $\varphi(M)$ cannot be contained
in any non-degenerate cone of $\mathbb{R}^n$. In particular, a
stochastically complete manifold cannot be minimally immersed into
a non-degenerate cone of $\mathbb{R}^n$.
\end{Corollario}

Note that even the statement of \cite{A} in its full generality
requires $|d\varphi|^2 \ge C >0$, an assumption that can be
overcome by the weighted requests \eqref{casoB}, \eqref{volumeB}.
Furthermore, in case $\beta=0$ we can replace $\Delta$-stochastic
completeness and the uniform control from below in \eqref{casoB}
with the properness of $\varphi$.

\begin{Corollario}\label{propria}
Let $(M,(,))$ be a Riemannian manifold. Then, there does not exist
any proper harmonic map $\varphi: M \rightarrow \mathbb{R}^n$,
such that $|d\varphi(x)|> 0$ on $M$ and $\varphi(M)$ is contained
into a non-degenerate cone of $\mathbb{R}^n$.
\end{Corollario}
\begin{Remark}
\emph{It is a well known open problem to deal with the case
$\theta= \pi/2$, that is, when the cone degenerates to a
half-space and the dimension $m$ is greater than $2$. When $m=2$,
$n=3$, by Hoffman-Meeks' half-space Theorem \cite{HM} the only
properly embedded minimal surfaces in a half-space are affine
planes. On the contrary, if $m\ge 3$ there exist properly embedded
minimal hypersurfaces even contained between two parallel
hyperplanes (the so called generalized catenoids). It is still an
open problem to find sufficient conditions on $M,\varphi$ in order
to have a Hoffman-Meeks' type result, and it seems quite difficult
to adapt the methods of the proof of \eqref{stimacono} for this
purpose.}
\end{Remark}

The next application of Theorem \ref{stimacono} has a topological
flavor. This result, which is interesting when $\varphi$ is not
proper, ensures that some kind of "patological" gathering around
points of $\overline{\varphi(M)}\backslash \varphi(M)$ does not
occur when the map is sufficiently well behaved. To make the
corollary more transparent, we state it using the sufficient
conditions \eqref{casoB} and \eqref{volumeB}. First we introduce
the following

\begin{Definizione}
Let $\mathcal{S}$ be a convex subset of $\mathbb{R}^n$. A point $p
\in \overline{\mathcal{S}}$ is called an $n$-corner of
$\mathcal{S}$ if it is the vertex of a non-degenerate cone
containing $\mathcal{S}$.
\end{Definizione}

\begin{Corollario}\label{inviluppo}
Let $(M,(,))$ be a complete Riemannian manifold and let $\varphi :
M \rightarrow \mathbb{R}^n$ be a map of class $C^2$. Suppose that
\eqref{casoB} holds, and that
\begin{equation}\label{tensione}
|\tau(\varphi)(x)| \le \frac{\widetilde{C}}{r(x)^\beta} \qquad \mbox{for }
r(x) \gg 1,
\end{equation}
for some $\widetilde{C}>0$ and $\beta \in \mathbb{R}$ as in
\eqref{casoB}. Assume also that \eqref{volumeB} holds. Then, the
convex envelope $\mathrm{Conv}(\varphi(M))$ contains no
$n$-corners.
\end{Corollario}
The second main theorem of this paper is obtained with a
modification of the idea used in the proof of Theorem
\ref{stimacono}. We obtain

\begin{Teorema}\label{stimacono2}
Let $\varphi : (M,(,)) \rightarrow \mathbb{R}^n$ be an isometric
immersion of an $m$-dimensional manifold satisfying the weak
maximum principle for the Hessian into a non-degenerate cone
$\mathcal{C}_{o,v,\theta}$. Assume the codimension restriction
\begin{equation}
0 < n-m < m
\end{equation}
and suppose that the sectional curvature of $M$ satisfy
\begin{equation}\label{sezionale}
^M\mathrm{Sect}_x \le \chi^2 \qquad \forall \ x \in M
\end{equation}
for some constant $\chi \ge 0$. Then
\begin{equation}\label{principale2}
\cos(\theta) \le \sqrt{d(\pi_v,\varphi(M)) \frac{\chi}{A_1}},
\end{equation}
where $A_1$ is as in \eqref{Aeta}.
\end{Teorema}
\begin{Remark}
\emph{With the notation $^M\mathrm{Sect} _x \le g(x)$ for some
function $g$ on $M$, we mean that, for every $x\in M$ and every
$2$-plane $\pi\le T_xM$, the sectional curvature of $\pi$ in $x$
satisfy $^M\mathrm{Sect}_x(\pi) \le g(x)$.}
\end{Remark}
As a consequence, we get the following corollaries: the former
generalizes results of Tompkins \cite{T}, Chern-Kuiper \cite{CK}
and Jorge-Koutroufiotis \cite{JK}, whereas the latter improves on
\cite{D}, Theorem $8.3$.

\begin{Corollario}\label{sezminzero}
Let $(M,\langle,\rangle)$ be a complete $m$-dimensional Riemannian
manifold with sectional curvature satisfying
\begin{equation}\label{sezionaleB}
- B^2 (1+r(x)^2) \prod_{j=1}^N \log^{(j)} r(x) \le \phantom{ } ^M\mathrm{Sect}_x
\le 0,
\end{equation}
for some $B>0$, $N \in \mathbb{N}$ and where $\log^{(j)}$ stands
for the $j$-iterated logarithm. Then, $M$ cannot be isometrically
immersed into a non-degenerate cone of $\mathbb{R}^{2m-1}$.
\end{Corollario}
\begin{Corollario}\label{kahler}
Let $(M^{2m},\langle,\rangle, J)$ be a K\"ahler manifold such that
the weak maximum principle for the Hessian holds. Then
$\varphi(M)$ cannot be isometrically immersed into a
non-degenerate cone of $\mathbb{R}^{3m-1}$.
\end{Corollario}

\begin{Remark}
\emph{The inequality on the left of \eqref{sezionaleB} provides a
sharp sufficient condition for the full Omori-Yau maximum
principle for the Hessian. Indeed, it would be sufficient that
this inequality hold along $2$-planes containing $\nabla r$ (see
Theorem 1.9 of \cite{PRS1}). As far as we know, it is an open
problem to obtain other general sufficient conditions ensuring the
validity of the weak maximum principle for the Hessian.}
\end{Remark}

\section{Proof of the Theorems}

\subsection{Proof of Theorem \ref{stimacono}}
First of all we observe that
$$
d(\pi_v,\varphi(M))= \inf_{x_o\in M} \langle \varphi(x_o),v\rangle,
$$
and that the right hand side of \eqref{principale} is invariant
under homothetic transformations of $\mathbb{R}^n$.
%
We choose $o$ as the origin of global coordinates, and for the
ease of notation we set
$$
b = \cos (\theta) \qquad  b \in (0,1).
$$
Furthermore, for future use, note that $\varphi(M)\subseteq
\mathcal{C}$ implies
\begin{equation}\label{defcono}
\langle \varphi(x),v\rangle \ge b|\varphi(x)| >0 \qquad \forall \
x \in M.
\end{equation}
Next, we reason by contradiction and we suppose that
\eqref{principale} is false. Therefore, there exists $x_o \in M$
such that
$$
\langle \varphi(x_o),v\rangle\sup_{x\in M}
\Big[\frac{|\tau(\varphi(x))|}{|d\varphi(x)|^2}\Big] < A_1b^2.
$$
By definition, and the fact that the inequality is strict, we can
find
$$
\xi \in (0,1) \quad , \quad \alpha \in \Big(0, \sqrt{(1-\xi)}\Big)
$$
such that
$$
\langle \varphi(x_o),v\rangle \sup_{x\in M}
\Big[\frac{|\tau(\varphi(x))|}{|d\varphi(x)|^2}\Big] < \Big(\xi
\alpha^2\sqrt{1-\alpha^2}\Big)b^2,
$$
thus:
\begin{equation}\label{primarichiestachi} \langle
\varphi(x_o),v\rangle |\tau(\varphi(x))| < \Big(\xi
\alpha^2\sqrt{1-\alpha^2}\Big)b^2|d\varphi(x)|^2 \qquad \forall \
x \in M.
\end{equation}
For the ease of notation we set $T=\langle
\varphi(x_o),v\rangle>0$ and $ a= b\alpha$; the last relation
becomes
\begin{equation}\label{richiestachi}
T|\tau(\varphi(x))| < \frac{\xi
a^2\sqrt{b^2-a^2}}{b}|d\varphi(x)|^2 \qquad \forall \ x \in M.
\end{equation}
Note also that
\begin{equation}\label{asopra}
a \in (0,b\sqrt{(1-\xi)}\big) \subseteq (0,b).
\end{equation}
Now, we define the following function:
$$
\displaystyle u(x) = \sqrt{T^2+a^2 |\varphi(x)|^2} -
\langle\varphi(x),v\rangle,
$$
and we note that, by construction, $u(x_o)> 0$. We first claim
that
\begin{equation}\label{limitato}
u < T \qquad \mathrm{on \ } M:
\end{equation}
indeed, an algebraic manipulation shows that \eqref{limitato} is
equivalent to
$$
\langle\varphi(x),v\rangle^2 +2T\langle\varphi(x),v\rangle - a^2
|\varphi(x)|^2  > 0 \qquad \mbox{on } M.
$$
On the other hand, using \eqref{defcono} the LHS of the above
inequality is bounded from below by $(b^2-a^2)|\varphi(x)|^2 >0$
since $a<b$, and the claim is proved.\par
We now consider the closed non empty set:
$$
\Omega_o = \{x\in M \; : \; u(x)\ge u(x_o) \}.
$$
Using \eqref{defcono} and the definition of $\Omega_o$ we deduce:
\begin{equation}\label{uno}
\sqrt{T^2+a^2 |\varphi(x)|^2} \ge b|\varphi(x)| + u(x_o).
\end{equation}
Since $u(x_o)>0$ by construction, we can square inequality
\eqref{uno} to obtain
\begin{equation}\label{ineqpolinom}
(b^2-a^2)|\varphi(x)|^2 + 2bu(x_o) |\varphi(x)| + u(x_o)^2-T^2 \le
0.
\end{equation}
Since $(b^2-a^2)>0$, the LHS of the above inequality is a
quadratic polynomial in $|\varphi(x)|$ with two distinct roots
$\alpha_- < 0 < \alpha_+$ (use Cartesio rule and
\eqref{limitato}), where $\alpha_\pm$ are given by
$$
\alpha_\pm = \big[b^2-a^2\big]^{-1} \Big\{ \pm \sqrt{(b^2-a^2)T^2
+ a^2u(x_o)^2}-bu(x_o)\Big\};
$$
therefore, \eqref{ineqpolinom} implies
%
%
%
%
\begin{equation}\label{sharp}
|\varphi(x)| \le \big[b^2-a^2\big]^{-1} \Big\{ \sqrt{(b^2-a^2)T^2
+ a^2u(x_o)^2}-bu(x_o)\Big\}
 \qquad \mbox{on  } \Omega_o.
\end{equation}
We then use the elementary inequality $\sqrt{1 +t^2} \le 1+t$ on
$\mathbb{R}_0^+$ to deduce
$$
\begin{array}{l}
\big[b^2-a^2\big]^{-1} \Big\{ \sqrt{(b^2-a^2)T^2 +
a^2u(x_o)^2}-bu(x_o)\Big\} \\[0.6cm]
=\displaystyle \frac{au(x_o)}{b^2-a^2} \sqrt{
1+\frac{(b^2-a^2)T^2}{a^2u(x_o)^2}} -\frac{bu(x_o)}{b^2-a^2}
\\[0.6cm]
\le \displaystyle \frac{au(x_o)}{b^2-a^2} \left(
1+\frac{T\sqrt{b^2-a^2}}{au(x_o)}\right) -\frac{bu(x_o)}{b^2-a^2}
\\[0.6cm]
= \displaystyle \frac{T}{\sqrt{b^2-a^2}} - \frac{u(x_o)}{b+a}
\end{array}
$$
and thus \eqref{sharp} together with $u(x_o)>0$ yields
\begin{equation}\label{upperbound}
|\varphi(x)| \le \frac{T}{\sqrt{b^2-a^2}} = \varphi_{max} \qquad
\mbox{on } \Omega_o.
\end{equation}
To compute $\Delta u$, we fix a local orthonormal frame $\{e_i\}$
and its dual coframe $\{\theta^i\}$. Then, writing
$du=u_i\theta^i$, a simple computation shows that
\begin{equation}
u_i = \frac{a^2\langle d\varphi(e_i),\varphi
\rangle}{\sqrt{T^2+a^2|\varphi|^2}} - \langle
d\varphi(e_i),v\rangle,
\end{equation}
and taking the covariant derivative we have $\nabla du =
u_{ij}\theta^i\otimes \theta^j$, where
$$
\begin{array}{lcl}
u_{ij} & = & \displaystyle - \frac{a^4\langle
d\varphi(e_i),\varphi \rangle \langle d\varphi(e_j),\varphi
\rangle}{\left(T^2+a^2|\varphi|^2\right)^{3/2}} - \langle \nabla d\varphi(e_i,e_j),v\rangle \\[0.7cm]
& & \displaystyle + \frac{a^2 \langle \nabla d\varphi(e_i,e_j),
\varphi \rangle + a^2\langle d\varphi(e_i),d\varphi(e_j)\rangle
}{\sqrt{T^2+a^2|\varphi|^2}}.
\end{array}
$$
Tracing the above expression we get
\begin{equation}\label{Deltau}
\Delta u = \langle \frac{S}{|\varphi|} \varphi - v, \tau(\varphi)
\rangle + S \frac{|d\varphi|^2}{|\varphi|} - \frac{1}{|\varphi|^2}
\frac{S^2}{\sqrt{T^2 + a^2|\varphi|^2}} \sum_{i=1}^m \langle
\varphi, d\varphi(e_i) \rangle^2
\end{equation}
on $M$, where we have defined
\begin{equation}\label{hacca}
S = S(x) = \frac{a^2 |\varphi(x)|}{\sqrt{T^2 +
a^2|\varphi(x)|^2}}.
\end{equation}
Note that, by \eqref{defcono},
\begin{equation}\label{Primo}
\Big| \frac{S}{|\varphi|} \varphi - v\Big|^2 \le S^2 -2bS + 1,
\end{equation}
and that
\begin{equation}\label{Secondo}
\sum_{i=1}^m \langle \varphi, d\varphi(e_i) \rangle^2 \le
\left\{\begin{array}{l} \displaystyle |\varphi|^2 \sum_{i=1}^m
|d\varphi(e_i)|^2 = |\varphi|^2 |d\varphi|^2; \\[0.3cm]
\displaystyle |\varphi|^2 = \frac{1}{m} |\varphi|^2 |d\varphi|^2
\qquad \mbox{if
 } \varphi \mbox{ is isometric.} \end{array} \right.
\end{equation}
The possibility, for the isometric case, of substituting $A_1$ with $A_m$ in \eqref{principale} depends only on
the above difference. Since the next passages are the same, we carry on with the general case. Substituting
\eqref{Primo}, \eqref{Secondo} in \eqref{Deltau} it follows that
\begin{equation}\label{ineqdelta}
\Delta u \ge - |\tau(\varphi)| \sqrt{S^2-2bS +1} + S
\frac{|d\varphi|^2}{|\varphi|} -
\frac{S^2}{\sqrt{T^2+a^2|\varphi|^2}} |d\varphi|^2.
\end{equation}
We now restrict our estimates on the RHS of \eqref{ineqdelta} on
$\Omega_o$. Then, \eqref{upperbound} holds and from
\eqref{richiestachi} we obtain
$$
\frac{|\tau(\varphi)|}{|d\varphi|^2} < \frac{\xi a^2
\sqrt{b^2-a^2}}{Tb} = \frac{\xi a^2}{\sqrt{T^2 + a^2
\varphi_{max}^2}} \le \frac{\xi a^2}{\sqrt{T^2 + a^2|\varphi|^2}}
= \frac{\xi S}{|\varphi|}.
$$
Inserting this inequality into \eqref{ineqdelta} we have
\begin{equation}\label{import}
\Delta u \ge \frac{a^2|d\varphi|^2}{\sqrt{T^2 +a^2|\varphi|^2}}
\Big[ 1 - \xi \sqrt{S^2-2bS+1} - \frac{S^2}{a^2} \Big].
\end{equation}
We want to find a strictly positive lower bound for $(1 - \xi
\sqrt{S^2-2bS+1} - S^2/a^2)$ on $\Omega_o$. Since $1-2bS+S^2$
represents a convex parabola and since $S$ is increasing in the
variable $|\varphi|$ on $[0,\varphi_{max}]$, its maximum is
attained either in $0$ or in $\varphi_{max}$. Since $S(0) = 0$,
$S(\varphi_{max}) = a^2/b>0$ we have
$$
S(\varphi_{max})^2 - 2b S(\varphi_{max}) + 1 = 1+
a^2\Big(\frac{a^2}{b^2} - 2\Big) < 1 = S(0)^2 - 2b S(0) + 1,
$$
thus we can roughly bound as follows:
$$
1 - \xi \sqrt{S^2-2bS+1} - \frac{S^2}{a^2} \ge 1 - \xi -
\frac{a^2}{b^2}
$$
and the RHS of the above inequality is strictly positive since $a
\in (0, b\sqrt{1-\xi})$. Therefore, \eqref{import} together with
\eqref{upperbound} yield
\begin{equation}\label{finale}
Lu = |d\varphi|^{-2} \Delta u \ge \frac{a^2}{\sqrt{T^2 +a^2|\varphi|^2}}\left[1-\xi-\frac{a^2}{b^2}\right]
\ge
\delta \qquad \mbox{on } \Omega_o,
\end{equation}
for some $\delta >0$.\\
There are now two possibilities:
\begin{itemize}
\item[$i)$] $x_o$ is an absolute maximum for $u$ on
    $M$. By assumption $|d\varphi(x_o)|^2>0$, and the finite
    form of the maximum principle yields $\Delta u(x_o) \le
    0$, so that $Lu(x_o)\le 0$. Since $x_o \in \Omega_o$
    \eqref{finale} immediately gives the contradiction.
\item[$ii)$] $\mathrm{Int}(\Omega_o) = \{ x \in M : u(x) >
    u(x_o) \} \neq \emptyset$. In this case, since
    $u(x)$ is bounded above on $M$, it is enough to evaluate
    inequality \eqref{finale} along a sequence $\{x_k\}$
    realizing the weak maximum principle for $L$, that is
    $u(x_k) > u^*-1/k$, $Lu(x_k) < 1/k$. Note that this
    sequence eventually lies in $\mathrm{Int}(\Omega_o)$.
\end{itemize}
This concludes the proof. $\blacktriangle$

%
%
%
%
%
\subsection{Proof of Corollary \ref{atsuji}}
If $M$ is stochastically complete and $|d\varphi|^2 \ge C$, then
it is straightforward to deduce that $M$ is $L$-stochastically
complete, where $L= |d\varphi|^{-2} \Delta $. Indeed, for every
$u\in C^2(M)$ with $u^*<+\infty$, along the sequence $\{x_k\}$
realizing the weak maximum principle for $\Delta$ we have also
$$
Lu(x_k) = |d\varphi(x_k)|^{-2}\Delta u(x_k) \le \frac{1}{Ck}.
$$
The result follows setting $\tau(\varphi) \equiv 0$ in Theorem
\ref{stimacono}. $\blacktriangle$

\subsection{Proof of Corollary \ref{propria}}
From \eqref{upperbound} in the proof of Theorem \ref{stimacono} we
deduce that $\varphi(\Omega_o)$ is bounded, hence
$\overline{\varphi(\Omega_o)}$ is compact. The properness
assumption implies that
$\varphi^{-1}(\overline{\varphi(\Omega_o)})$ is compact, thus
$\Omega_o$ is compact. Therefore, it is enough to use the finite
form of the maximum principle in \eqref{finale}. $\blacktriangle$
\subsection{Proof of Corollary \ref{inviluppo}}
We reason by contradiction and let $p \in
\mathrm{Conv}(\varphi(M))$ be an $n$-corner.
\begin{itemize}
\item[-] If $p\in \mathrm{Conv}(\varphi(M)) \backslash
\overline{\varphi(M)}$ fix a small ball around $p$ contained in
$\mathbb{R}^n\backslash\overline{\varphi(M)}$, and cut the corner
transversally with an hyperplane sufficiently near to $p$; it is
immediate to see that in this way we produce a convex set
containing $\varphi(M)$ and strictly smaller than
$\mathrm{Conv}(\varphi(M))$, contradiction.
\item[-] Suppose now $p \in \varphi(M)$, and let $x\in M$ such
that $\varphi(x)=p$. Consider the map $d\varphi_{|x}$; by
assumption, there exists a direction $v \in T_xM$ such that
$|d\varphi_{|x} v| \neq 0$, thus by continuity we can take a curve
$$
\gamma : (-\varepsilon, \varepsilon) \rightarrow M \quad , \quad
\gamma(0)=x \quad , \quad \gamma'(0) = v
$$
with $\varepsilon$ small such that
$|d\varphi_{|\gamma(t)}(\gamma'(t))| \neq 0$ on $(-\varepsilon,
\varepsilon)$. Therefore, $\varphi\circ \gamma$ is an immersed
curve in $\mathbb{R}^n$, and this fact
contradicts the assumption that $p$ is an $n$-corner.\\
\item[-] If $p \in \overline{\varphi(M)}\backslash
    \varphi(M)$, choose $\pi_v$ as the hyperplane orthogonal
    to the direction of the cone and passing through $p$. It
    follows that $d(\varphi(M), \pi_v)=0$. By \eqref{casoB}
    and \eqref{tensione}, we argue that
    $|\tau(\varphi)|/|d\varphi|^2$ is bounded above on $M$. By
    Remark \ref{sufficienticond}, \eqref{casoB} and
    \eqref{volumeB} ensure that $M$ is $L$-stochastically
    complete, where $L=|d\varphi|^{-2}\Delta$. By Theorem
    \ref{stimacono} we conclude the validity of
    \eqref{principale} which gives $\theta= \pi/2$,
    contradiction. $\blacktriangle$\\
\end{itemize}
\subsection{Proof of Theorem \ref{stimacono2}}
We follow the proof of Theorem \ref{stimacono} verbatim
substituting \eqref{primarichiestachi} with
\begin{equation}\label{primarichiestachiB}
\langle \varphi(x_o), v \rangle \chi < \Big(\xi \alpha^2 \sqrt{1-\alpha^2}\Big) b^2
\end{equation}
for some $\xi \in (0,1), \alpha \in (0, \sqrt{1-\xi})$, \eqref{richiestachi} with
\begin{equation}\label{richiestachiB}
T \chi < \frac{\xi a^2 \sqrt{b^2-a^2}}{b}
\end{equation}
and arrive up to inequality \eqref{upperbound} included. Next, we
fix $x \in \Omega_o$ and we let $X,Y \in T_xM$  be orthonormal
vectors spanning the $2$-plane $\pi$. From Gauss equations and
\eqref{sezionale} we have
\begin{equation}\label{gauss}
\langle II(X,X) ,  II(Y,Y) \rangle - | II(X,Y)|^2 =
^M\mathrm{Sect}(\pi) \le \chi^2,
\end{equation}
where $II$ is the second fundamental form. Since $0<n-m<m$, by
Otsuki lemma (see \cite{D}, Lemma $3.1$) it follows that there
exists a unit vector $W \in T_xM$ such that
$$
|II(W,W)| \le \chi
$$
hence, from \eqref{richiestachiB} and \eqref{upperbound} we deduce
\begin{equation}\label{otsuki}
|II(W,W)| < \frac{\xi a^2 \sqrt{b^2-a^2}}{Tb} = \frac{\xi
a^2}{\sqrt{T^2+a^2\varphi_{max}^2}} \le \frac{\xi
a^2}{\sqrt{T^2+a^2|\varphi(x)|^2}}.
\end{equation}
Next, we let $\gamma: [0,\varepsilon) \rightarrow M$, $\varepsilon >0$, be the geodesic characterized by the
initial data
$$
\gamma(0) = x \quad , \quad \gamma'(0) = W.
$$
Call $s \in [0,\varepsilon)$ the arc-length parameter and define
the function
$$
g : [0,\varepsilon) \rightarrow \mathbb{R} \qquad g(s) = u
(\gamma(s)).
$$
A simple computation, using the fact that $\varphi$ is an
isometric immersion, gives:
\begin{equation}\label{gsec}
g''(s) = \langle \frac{S}{|\varphi(\gamma)|}\varphi(\gamma) - v,
II(\gamma',\gamma') \rangle + \frac{S}{|\varphi(\gamma)|} -
\frac{S^3}{a^2 |\varphi(\gamma)|^3} \langle d\varphi(\gamma'),
\varphi(\gamma) \rangle^2,
\end{equation}
where $S$ has the expression in \eqref{hacca}, evaluated at $x=
\gamma(s)$. Since
$$
\Big| \frac{S}{|\varphi|}\varphi - v\Big|^2 \le 1+S^2-2bS \quad ,
\quad \langle d\varphi(\gamma'), \varphi \rangle^2 \le
|d\varphi(\gamma')|^2|\varphi|^2 = |\varphi|^2
$$
Setting $S_o = S(\gamma(0))$, evaluating at $s=0$ we deduce
\begin{equation}\label{semifinale}
g''(0) \ge - |II(W,W)| \sqrt{S_o^2-2bS_o + 1} + \frac{a^2 -
S_o^2}{\sqrt{T^2+a^2|\varphi(\gamma)|^2}}.
\end{equation}
Inserting \eqref{otsuki} into \eqref{semifinale} we get
\begin{equation}
g''(0) \ge \frac{a^2}{\sqrt{T^2+a^2|\varphi(\gamma)|^2}} \Big[ 1 -
\xi \sqrt{1+S_o^2-2bS_o} - \frac{S_o^2}{a^2}\Big].
\end{equation}
Proceeding as in the proof of Theorem \ref{stimacono}, since $a
\in (0,b\sqrt{1-\xi}) \subset (0,b)$
$$
g''(0) \ge \frac{a^2}{\sqrt{T^2+a^2|\varphi(\gamma)|^2}} \Big[ 1-\xi - \frac{a^2}{b^2}\Big] \ge
\frac{a^2\sqrt{b^2-a^2}}{bT}[1-\xi-\frac{a^2}{b^2}] = \delta > 0,
$$
where $\delta$ is independent from $x\in \Omega_o$ and from
$W$.\\
On the other hand, a standard computation using the fact that
$\gamma$ is a geodesic and the definition of the Hessian of a
function, gives $g''(0) = \mathrm{Hess}_x u(W,W)$. Putting
together the last two inequalities we obtain
\begin{equation}\label{contraddiz}
\mathrm{Hess}_x u(W,W) \ge \delta > 0.
\end{equation}
If $x_o$ is an absolute maximum of $u$, then from
\eqref{contraddiz} we immediately contradict the finite maximum
principle, otherwise
\begin{equation}
\mathrm{Int}(\Omega_o) = \{ x \in M : u(x) > u(x_o)\} \neq
\emptyset
\end{equation}
and \eqref{contraddiz} gives
\begin{equation}
\inf_{x \in \mathrm{Int}(\Omega_o)} \sup_{\scriptsize {\begin{array}{c} Y \in T_x M \\
|Y| = 1 \end{array}}} \mathrm{Hess}_x u (Y,Y) \ge \delta > 0,
\end{equation}
contradicting the validity of the weak maximum principle for the
Hessian operator since the function $u$ in bounded above on $M$.
This completes the proof of Theorem \ref{stimacono2}.
$\blacktriangle$
\subsection{Proof of Corollary \ref{sezminzero}}
By \eqref{sezionaleB}, using Theorem $1.9$ of \cite{PRS1} we have the validity of the weak maximum principle for
the Hessian. The result follows immediately setting $\chi=0$ and $n=2m-1$. $\blacktriangle$
\subsection{Proof of Corollary \ref{kahler}}
The proof follows the same lines as in \cite{D}, so we only sketch
it. From the assumptions, since the codimension is $m-1 < m$, for
every $x \in M$ the theory of flat bilinear forms ensure the
existence of a vector $Z \in T_xM$, $|Z|=1$ such that $II(JZ,JZ) =
- II(Z,Z)$. We define $u, \Omega_o$ as in Theorem
\ref{stimacono2}. Expression \eqref{gsec} gives at every point
$x$, and for every $X \in T_xM$, $|X|=1$
$$
\mathrm{Hess}_x u(X,X) \ge \langle
\frac{S}{|\varphi(x)|}\varphi(x) - v, II(X,X) \rangle +
\frac{S}{|\varphi(x)|} \Big( 1 - \frac{S^2}{a^2} \Big).
$$
This calculation is independent from the value of $a \in (0,b)$.
If $a$ is chosen to be sufficiently small that $S^2/a^2 < \delta <
1$ (note that, by definition, on $\Omega_o$ it holds $S=O(a^2)$
and $S/|\varphi| \ge a^2/T$), evaluating along a sequence
$\{x_k\}$ satisfying the weak maximum principle for the Hessian we
deduce, for $k$ sufficiently large,
$$
\begin{array}{lcl}
\displaystyle \langle \frac{S}{|\varphi(x_k)|}\varphi(x_k) - v,
II(X_k,X_k) \rangle & \le & \displaystyle \mathrm{Hess}_x
u(X_k,X_k) - \frac{S}{|\varphi(x_k)|} (1-\delta) \\[0.2cm]
& \le & \displaystyle \frac{1}{k} - \frac{a^2}{T} (1-\delta) < 0

\end{array}
$$
for every $X_k \in T_{x_k}M, \ |X_k|=1$. This fact contradicts the existence of $Z$. $\blacktriangle$
\section{Ending remarks}
We end the paper with some observations regading Theorem
\ref{stimacono2} and Corollary \ref{sezminzero}. In particular, we
stress the difference of applicability between the weak maximum
principles for the Laplacian and for the full Hessian operator. A
first striking difference is pointed out by Proposition 40 of
\cite{PRS4}, which states that every Riemannian manifold
satisfying the weak maximum principle for the Hessian must be
non-extendible (that is, non isometric to any proper open subset
of another Riemannian manifold). For example, for every Riemannian
manifold $M$ and $p\in M$, $M\backslash\{p\}$ does not satisfy the
weak maximum principle for the Hessian. This is in sharp contrast
to Theorem IX.$3$ of \cite{C}, whose immediate application is the
following
\begin{Proposizione}
Let $(M,\langle,\rangle)$ be a stochastically complete
$m$-dimensional manifold with Riemannian measure $\mu$, and let
$K\subset M$ be a compact submanifold of dimension $k\ge 0$ such
that $M\backslash K$ is connected. Then, if $m-k \ge 2$,
$M\backslash K$ is stochastically complete.
\end{Proposizione}
\noindent \textbf{Proof:}\\
Denoting with $B_\varepsilon= \{x : d(x,K)<\varepsilon\}$ and with
$p_\varepsilon$ the Dirichlet heat kernel on $M\backslash
B_\varepsilon$, by Theorem IX.$1$ of \cite{C} $p_\varepsilon
\uparrow p$ uniformly on compact subsets of $M\backslash K \times
[0,+\infty)$. Therefore, the Dirichlet heat kernel $\widetilde{p}$
of $M\backslash K$ coincides with $p_{|M\backslash K}$, and since
$\mu(K)=0$ the Brownian motion $X_t$ on $M\backslash K$ satisfies
$$
\mathbb{P}(X_t \in M\backslash K \ | \ X_0=x) = \int_{M\backslash
K} \widetilde{p}(x,y,t)d\mu(y) = \int_M p(x,y,t)d\mu(y) = 1
$$
for every $t>0$, $x\in M\backslash K$, and this shows that
$M\backslash K$ is stochastically complete. $\blacktriangle$\\
\\
For those familiar with stochastic
calculus and potential theory this proposition may look almost
trivial since it turns out that, in these assumptions,
the set $K$ is polar with respect to the Brownian motion (i.e. it has zero capacity).\\
Since, by Theorem 1.9 of \cite{PRS1}, geodesic completeness and a
well-behaved sectional curvature imply the full Omori-Yau maximum
principle for the Hessian, one might ask if, keeping a
well-behaved sectional curvature and relaxing geodesic
completeness to the property of non-extendibility, one could prove
the validity of the weak maximum principle for the Hessian. This
is false, as the following simple counterexample shows. Consider
on Euclidean space $\mathbb{R}^3$ the standard cone
$$
M = \{ x=(x_1,x_2,x_3) \neq (0,0,0) : x_3=\sqrt{x_1^2+x_2^2} \}.
$$
In polar coordinates $(r,\theta)$, where $r\in (0,+\infty)$ and
$\theta\in [0,2\pi)$, the cone can be parametrized as $x_1 =
r\cos\theta$, $x_2 = r\sin\theta$, $x_3=r$. Therefore, the induced
metric reads
$$
ds^2 = 2dr^2 + r^2d\theta^2;
$$
this shows that the cone is trivially non-extendible as a
Riemannian manifold (every such extension $N$ must contain only
one point not in $M$, but the metric is singular in $r=0$) but,
since $M$ is a flat embedded hypersurface trivially contained into
a non-degenerate cone, from Theorem \ref{stimacono2} the weak
maximum principle for the Hessian necessary fails. We conclude the
remark observing that $M$ is stochastically complete. Indeed, from
the form of the metric we deduce that the normal projection onto
the hyperplane $x_3=0$ gives a quasi-isometry between $M$ and
$\mathbb{R}^2\backslash\{(0,0)\}$, preserving divergent sequences
and such that the derivatives of the metric on $M$ are controlled
by those of $\mathbb{R}^2\backslash \{(0,0)\}$. Therefore,
stochastic completeness follows applying a slight modification of
Proposition $3.4$ in \cite{PRS1}.\\
\vspace{0.5cm}

\textbf{Added in proof.} With the article in press, we found that
the very recent paper \cite{RM} is deeply related both to our work
and to Atsuji's paper \cite{A}. Indeed, the author has succeeded
in recovering Atsuji's result with the aid of a geometrical
approach via the weak maximum principle very similar to the one
presented here and in \cite{PRS1}.


\begin{thebibliography}{biblio}

\bibitem[A]{A} A. Atsuji, \emph{Remarks on harmonic maps into a
cone from a stochastically complete manifold.} Proc. Japan Acad.
\textbf{75} Ser. A (1999), 105-108.

\bibitem[AG]{AG} L.J. Al\'ias, S.C. Garc\'ia-Mart\'inez, \emph{On the
scalar curvature of constant mean curvature hypersurfaces in space
forms.} J. Math. Anal. Appl. \textbf{363} (2010), no. 2, 589-597.

\bibitem[APD]{APD} L.J. Al\'ias, G. Pacelli Bessa, M. Dajczer,
\emph{The mean curvature of cylindrically bounded submanifolds.}
Math. Ann. \textbf{345} (2009), no. 2, 367-376.

\bibitem[BK]{BK} C. Baikoussis, T. Kouforgiorgos, \emph{Harmonic maps
into a cone.} Arch. Math. \textbf{40} (1983), 372-386.

\bibitem[C]{C} I. Chavel, \emph{Eigenvalues in Riemannian
Geometry.} Academic Press, Orlando, FL, 1984.

\bibitem[CK]{CK} S.S. Chern, N.H. Kuiper, \emph{Some theorems on
the isometric imbedding of compact Riemannian manifolds in
Euclidean space.} Annals of Math. \textbf{56} (1952), 422-430.

\bibitem[D]{D} M. Dajczer, \emph{Submanifolds and Isometric immersions.}
Math. Lecture Series vol \textbf{13}, Publish or Perish, Houston,
TX, 1990.

\bibitem[G]{G} A. Grigor'yan, \emph{Analytic and geometric background of
recurrence and non-explosion of the Brownian motion on Riemannian
manifolds.} Bull. Am. Math. Soc. \textbf{36} (1999), 135-249.

\bibitem[GW]{GW} R. Greene, H.H. Wu, \emph{Function theory on manifolds
which possess a pole.} LNM \textbf{699}, Springer-Verlag, Berlin
1979.

\bibitem[HM]{HM} D. Hoffman, W. Meeks, \emph{The strong half-space
Theorem for minimal surfaces.} Invent. Math. \textbf{101} (1990),
373-377.

\bibitem[JK]{JK} L. Jorge, D. Koutroufiotis, \emph{An estimate for
the curvature of bounded submanifolds.} Amer. J. of Math.
\textbf{103} (1981), 711-725.

\bibitem[O]{O} H. Omori, \emph{Isometric immersions of Riemannian
manifolds.} J. Math. Soc. Japan \textbf{19} (1967), 205-214.

\bibitem[PRS1]{PRS1} S. Pigola, M. Rigoli, A.G. Setti, \emph{Maximum
principles on Riemannian manifolds and applications.} Memoirs of
the Am. Math. Soc. vol.\textbf{174} (2005).

\bibitem[PRS2]{PRS2} S. Pigola, M. Rigoli, A.G. Setti, \emph{A remark on
the maximum principle and stochastic completeness.} Proc. Am.
Math. Soc. \textbf{131} (2003), 1283-1288.

\bibitem[PRS3]{PRS3} S. Pigola, M. Rigoli, A.G. Setti, \emph{Volume
growth, "a priori" estimates and geometric applications.} Geom.
Funct. An. \textbf{13} (2003), 1302-1328.

\bibitem[PRS4]{PRS4} S. Pigola, M. Rigoli, A.G. Setti, \emph{Maximum
principles at infinity on Riemannian manifolds: an overview.}
Matem\'atica Contempor\^anea \textbf{31} (2006), 81-128.

\bibitem[RM]{RM} A. Ranjbar-Motlagh, \emph{On harmonic maps from stochastically
complete manifolds.} Arch. Math. (Basel) \textbf{92} (2009), no.
6, 637-644.

\bibitem[RSV]{RSV} M. Rigoli, M. Salvatori, M. Vignati, \emph{Some
remarks on the weak maximum principle.} Rev. Mat. Iberoam.
\textbf{21} (2005), 459-481.

\bibitem[S]{S} K.T. Sturm, \emph{Analysis on local Dirichlet spaces I.
Recurrence, conservativeness and $L^p$-Liouville properties.} J.
Reine Angew. Math. \textbf{456} (1994), 173-196.

\bibitem[T]{T} C. Tompkins, \emph{Isometric embedding of flat
manifolds in Euclidean space.} Duke Math. J. \textbf{5} (1939),
58-61.

\end{thebibliography}
\end{document}